\documentclass[pdflatex,sn-basic]{sn-jnl}

\usepackage{graphicx}%
\usepackage{multirow}%
\usepackage{amsmath,amssymb,amsfonts}%
\usepackage{amsthm}%
\usepackage[title]{appendix}%
\usepackage{xcolor}%
\usepackage{textcomp}%
\usepackage{manyfoot}%
\usepackage{booktabs}%
\usepackage{algorithm}%
\usepackage{algorithmicx}%
\usepackage{algpseudocode}%
\usepackage{listings}%

\theoremstyle{thmstyleone}%
\newtheorem{theorem}{Theorem}
\newtheorem{proposition}{Proposition}%
\newtheorem{lemma}{Lemma}
\newtheorem{corollary}{Corollary}

\theoremstyle{thmstyletwo}%
\newtheorem{remark}{Remark}%

\theoremstyle{thmstylethree}%
\newtheorem{definition}{Definition}%

\usepackage[T1]{fontenc}

\newcommand{\EE}{\mathbb{E}}
\newcommand{\VV}{\mathbb{V}}
\newcommand{\CC}{\text{Cov}}

\newcommand{\ba}{\boldsymbol{\alpha}}
\makeatletter
\def\widebreve{\mathpalette\wide@breve}
\def\wide@breve#1#2{\sbox\z@{$#1#2$}%
     \vbox{\m@th\ialign{##\crcr
\kern0.08em\brevefill#1{0.8\wd\z@}\crcr\noalign{\nointerlineskip}%
                    $\hss#1#2\hss$\crcr}}}
\def\brevefill#1#2{$\m@th\sbox\tw@{$#1($}%
  \hss\resizebox{#2}{\wd\tw@}{\rotatebox[origin=c]{90}{\upshape(}}\hss$}
\makeatother

\begin{document}

\title[Asymptotic efficiency comparison of \\explicit Dirichlet shape estimators]{Asymptotic efficiency comparison of \\explicit Dirichlet shape estimators}


\author*[1]{\fnm{Ioannis} \sur{Oikonomidis}}\email{goikon@math.uoa.gr}

\author[1]{\fnm{Samis} \sur{Trevezas}}\email{strevezas@math.uoa.gr}

\affil*[1]{\orgdiv{Department of Mathematics}, \orgname{National and Kapodistrian University of Athens}, \orgaddress{\street{Panepistimioupolis}, \city{Athens}, \postcode{15784}, \country{Greece}}}


\abstract{In this study, we compare the asymptotic efficiency of three closed-form estimators of the Dirichlet shape parameters. We start by proving that, in the beta distribution, the score-adjusted moment estimator (SA) has strictly smaller componentwise asymptotic variances than the usual moment estimator (M1) for both parameters. We then introduce a new moment estimator (M2) that combines the marginal moment equations to obtain a common estimate of the total concentration. M2 coincides with M1 in the beta case and has strictly smaller componentwise asymptotic variances for every positive shape vector in three or more dimensions. We then show that SA is asymptotically more efficient than M2 when all shape parameters are equal; if at least two parameters differ, each estimator has smaller asymptotic variance than the other for some linear combination of the parameters. Furthermore, SA has strictly smaller componentwise asymptotic variances than M2 when all shape parameters are at least 1. Finally, we show that this ordering can reverse if a parameter becomes too small, and obtain relevant mathematical and computational bounds. Overall, the study shows that, for the Dirichlet distribution and under a broad range of parameter settings, score-adjusted estimation outperforms the classic moment estimation and identifies the limits of that advantage.}


\maketitle

\section{Introduction}

The Dirichlet distribution is a multivariate extension of the beta distribution from $[0, 1]$ to the standard simplex \citep{kotz2000continuous}. It is therefore a natural model for compositional data and category probabilities \citep{pal2025revisiting, bilancia2025stochastic}. The maximum likelihood estimator (MLE) of its shape parameter vector is not available in closed form and requires the use of Newton--Raphson or related numerical methods \citep{wicker2008maximum}. Explicit estimators are therefore useful either as initial values for iterative likelihood procedures or as computationally inexpensive alternatives when estimation must be performed repeatedly \citep{kim2026novel}. This is an active research topic for other common distribution families as well, such as the Weibull \citep{jokiel2024estimation} and the negative binomial \citep{zhao2023new}. For the Dirichlet, explicit likelihood-based estimation methods have also been proposed. Specifically, \citet{wicker2008maximum} worked on an explicit likelihood approximation for initializing numerical optimization and \citet{chang2023asymptotically} applied a likelihood-based correction to a moment estimator. 

In this study, we focus on closed-form moment-type estimators and compare their asymptotic efficiency, continuing our previous work \citep{oikonomidis2025moment}. \citet{tamae2020score} introduced the score-adjusted moment estimator (SA) for the beta distribution and proved that, when the two shape parameters are equal, SA has smaller asymptotic variances than the usual moment estimator (M1) for both parameters. Their numerical comparisons suggested that the inequality also holds for unequal parameter values. In \citet{oikonomidis2025moment}, we generalized SA to the Dirichlet distribution and proved that, as the dimension tends to infinity, the score-adjusted estimator has smaller componentwise asymptotic variances than the common moment estimator. In this study, we keep the dimension fixed and show that SA uniformly dominates M1 componentwise in the beta distribution, extending the result of \citet{tamae2020score}. 

M1 uses the first two marginal moments of each component separately \citep{minka2000estimating, oikonomidis2025moment}. We introduce a new moment estimator (M2) that uses all components to estimate the total concentration. A similar moment estimator was considered by \citet{ng2011dirichlet}, using all but one of the components. We show that M1 and M2 coincide for the beta distribution and that M2 has strictly smaller asymptotic variance than M1 for every shape parameter in three or more dimensions. We therefore use M2 as the basis for comparison with SA.

SA and M2 both estimate the mean composition by the sample mean vector, but they estimate the total concentration differently. In this study, we prove that SA is asymptotically more efficient than M2 when all shape parameters are equal, and that the difference of the asymptotic variance-covariance matrices is indefinite in any other case, extending the result of \citet{tamae2020score}. We therefore focus on the componentwise asymptotic variances and prove that if all parameters are at least $1$, then SA outperforms M2. However, this ordering is not uniform, and simple counterexamples can be found. We then move on to derive both analytical and computational bounds for the infimum value the parameters can attain for SA to have smaller componentwise asymptotic variances than M2. Overall, the study shows that under a broad range of parameter settings, SA estimation outperforms the classic moment estimation and identifies the limits of that advantage.

The remainder of the paper is organized as follows. Section~\ref{sec:background} reviews the necessary background on moment-type estimators of the Dirichlet shape parameters. Section~\ref{sec:SA-M1} focuses on the beta distribution and proves that SA uniformly outperforms the classical moment estimator M1. Section~\ref{sec:M2} introduces the new moment estimator M2, derives its asymptotic covariance matrix, and compares it with M1. Section~\ref{sec:SA-M2} compares SA with M2 and derives sufficient conditions under which SA outperforms M2. Section~\ref{sec:numerical} presents theoretical and computational analyses of the boundaries associated with these sufficient conditions. Finally, Section~\ref{sec:conclusion} provides concluding remarks and directions for future research.

\section{Background}\label{sec:background}
\noindent
In this section, we give a brief overview of the Dirichlet distribution and the moment-type estimators of its parameters. The results are stated without proof, as our aim is to provide the necessary background to understand the findings of this study. A more thorough treatment is available in \citet{oikonomidis2025moment}. 

\begin{definition}\label{def:dirichlet}
    A random vector $\mathbf{X}=(X_1,\ldots, X_k) \in \mathbb{R}^{k}$, $k\geq 2$, is said to follow the Dirichlet distribution with shape parameters $\boldsymbol{\alpha} \in (0, +\infty)^k$, denoted by $\mathbf{X}\sim\mathcal{D}_k(\boldsymbol{\alpha})$, if under the parameterization $x_k=1-\sum_{i=1}^{k-1}x_i$ it has a probability density function with respect to the Lebesgue measure $dx_1\cdots dx_{k-1}$ on $\mathbb{R}^{k-1}$ of the form
    \begin{equation*}
        f(\mathbf{x};\boldsymbol{\alpha}) = \frac{\Gamma(\alpha_0)}{\prod_{i=1}^k\Gamma(\alpha_i)}\,\prod_{i=1}^k x_i^{\alpha_i-1} \, \, \mathbb{I}_{\mathbf{S}^{\circ}_{k-1}}(\mathbf{x}),
    \end{equation*}
    where $\mathbf{S}^{\circ}_{k-1}$ is the open $(k-1)$-dimensional simplex, i.e., $\mathbf{S}^{\circ}_{k-1} := \{\mathbf{x} \in (0,1)^k: \sum_{i=1}^k x_i = 1\}$, $\Gamma$ is the gamma function, and $\alpha_0 := \sum_{i=1}^k \alpha_i$.
\end{definition}
\noindent
To keep the notation concise, we adopt a number of conventions. First, we will omit the parameter vector $\boldsymbol{\alpha}$ from the moment notation, thus writing $\EE(\mathbf{X})$ instead of $\EE_{\boldsymbol{\alpha}}(\mathbf{X})$. To highlight the fact that expectation and variance are functions of the parameters, we will instead write $e_{h}(\boldsymbol{\alpha}):=\EE_{\boldsymbol{\alpha}}[h(\mathbf{X})]$ and $v_{h}(\boldsymbol{\alpha}):=\VV_{\boldsymbol{\alpha}}[h(\mathbf{X})]$ whenever necessary. Furthermore, the observations of a random sample of size $N$ will be indexed by $n \in [N]:=\{1, 2, \dots, N\}$, and the $k$ elements of a vector will be indexed by $i \in [k]$. Therefore, we will denote the mean of a random sample $\mathbf{X}_n = (X_{1n}, \dots, X_{kn}), \, n \in [N]$, by $\overline{\mathbf{X}}$, and its components by $\overline{X_i}, \, i \in [k]$. 

We now give the definition of the general class of moment-type estimators.

\begin{definition}\label{def:me}
    Let $\mathbf{X}_n = (X_{1n}, \dots, X_{kn}), \, n \in [N]$, be a random sample from $\mathcal{D}_k(\boldsymbol{\alpha})$. Let $h:\mathbf{S}^{\circ}_{k-1}\rightarrow\mathbb{R}^m$ be a function such that $e_h: \mathbb{R}_+^k \rightarrow\mathbb{R}^m$ is an injection, and $g:\mathbb{R}^m\rightarrow \mathbb{R}^k$ be a function such that $g|_{e_h\left(\mathbb{R}_+^k\right)} = e_h^{-1}$. Then, the estimator
    \[
        \widetilde{\boldsymbol{\alpha}}_N = g\left[\overline{h\left(\mathbf{X}\right)}\right],
    \]
    is called a moment-type estimator of $\boldsymbol{\alpha}$. The estimator exists if $\widetilde{\boldsymbol{\alpha}}_N\in \mathbb{R}^k_{+}$.
\end{definition}
\begin{proposition} \label{prop:me-asymp}
    If $g$ is continuously differentiable in a neighborhood of $\boldsymbol{\mu}:=e_h(\boldsymbol{\alpha})$ and $\EE\left(||h(\mathbf{X})||^2\right)$ is finite, then the moment-type estimator $\widetilde{\boldsymbol{\alpha}}_N$ exists with probability tending to 1, and
    \[
        \sqrt{N}\left(\,\widetilde{\boldsymbol{\alpha}}_N-\boldsymbol{\alpha}\,\right) \overset{\mathcal{L}}{\longrightarrow} \mathcal{N}_k\left(\mathbf{0}, \boldsymbol{\Sigma}\right),
    \]
    where $\boldsymbol{\Sigma} := \mathbf{G}\mathbf{V}\mathbf{G}^\top$, $\mathbf{G} :=\nabla g(\boldsymbol{\mu})$, and $\mathbf{V} := v_h(\boldsymbol{\alpha})$.
\end{proposition}
As stated in Proposition \ref{prop:me-asymp}, these estimators are consistent and asymptotically normal. We therefore compare their large-sample performance in terms of relative asymptotic efficiency. Even though the MLE is asymptotically efficient---and, in the Dirichlet case, a moment-type estimator \citep{oikonomidis2025moment}---it is not available in closed form. Specifically, for a random sample $\mathbf{X}_n,\, n \in [N]$, from $\mathcal{D}_k(\boldsymbol{\alpha})$, the score equation system takes the form: 
\begin{equation*}
    \psi(\alpha_i) - \psi\left(\alpha_0\right) = \frac{1}{N} \sum_{n=1}^N \ln X_{in}, \quad i = 1, \ldots, k,
\end{equation*}
where $\psi$ is the digamma function. 

A multidimensional estimator with asymptotic variance--covariance matrix $\boldsymbol{W}$ is asymptotically more efficient than another one with respective matrix $\boldsymbol{\Sigma}$, if the difference $\boldsymbol{\Sigma}-\boldsymbol{W}$ is positive semidefinite; equivalently, the more efficient estimator has equal or smaller asymptotic variance for any linear combination of the estimated parameters. Nevertheless, such an ordering need not hold. Therefore, we also consider the weaker notion of componentwise efficiency, which compares only the corresponding diagonal elements. Under this criterion, one estimator is characterized as \textit{componentwise asymptotically more efficient} than another one, if it has an equal or smaller asymptotic variance for every estimated component. 

We study three closed-form moment-type estimators. The moment estimator M1 and the score-adjusted moment estimator SA are presented in Propositions \ref{prop:dirichlet-me1} and \ref{prop:dirichlet-sa}, respectively \citep{oikonomidis2025moment}. We introduce a second moment estimator (M2) in Section \ref{sec:M2}. The three estimators are well-defined and positive almost surely for $N\geq2$.

\begin{proposition} \label{prop:dirichlet-me1}
    Let $\mathbf{X}_n, \, n \in [N], \,N\geq2$, be a random sample from $\mathcal{D}_k(\boldsymbol{\alpha})$. The Moment Estimator
    \begin{equation*}
        \widetilde{\alpha}_{N,i} = \frac{\overline{X_i} \left(\overline{X_i} - \overline{X_i^2}\right)}{\overline{X_i^2}-\overline{X_i}^2}, \qquad i \in [k], 
    \end{equation*}
    is a moment-type estimator with asymptotic variance-covariance matrix $\boldsymbol{\Sigma}$ such that
    \begin{equation}\label{eq:Sij}
    \boldsymbol{\Sigma}_{ij} = \frac{(\alpha_i + 1)(\alpha_j + 1)(\delta_{ij}\alpha_0-\alpha_i)\alpha_j\alpha_0}{(\alpha_0-\alpha_i)(\alpha_0-\alpha_j)(\alpha_0+2)(\alpha_0+3)} \left\{ \frac{2(\alpha_0+1)^2}{\alpha_j+1}\delta_{ij}-\frac{2\alpha_0^2+\alpha_0+1}{\alpha_0+1} \right\},
    \end{equation}
    for $i,j \in [k]$, where $\delta_{ij}$ is the Kronecker delta. Its $i$-th diagonal asymptotic variance is
    \begin{equation}\label{eq:M1-diag}
    \sigma_{i}^2(\boldsymbol{\alpha})
    =
    \frac{\alpha_0\alpha_i(\alpha_i+1)
    \left\{2(\alpha_0+1)^3-(\alpha_i+1)(2\alpha_0^2+\alpha_0+1)\right\}}
    {(\alpha_0+1)(\alpha_0+2)(\alpha_0+3)(\alpha_0-\alpha_i)}, \quad i \in [k].
    \end{equation}
\end{proposition}

\begin{proposition} \label{prop:dirichlet-sa}
    Let $\mathbf{X}_n, \, n \in [N], \,N\geq2$, be a random sample from $\mathcal{D}_k(\boldsymbol{\alpha})$. The Score-Adjusted Moment Estimator
    \begin{equation*}
        \widebreve{\alpha}_{N,i} = \frac{(k-1)\overline{X_i}}{\sum_{m=1}^k\left(\,\overline{X_m \ln X_m}-\overline{X_m} \, \overline{\ln X_m}\,\right)}, \quad i \in [k],
    \end{equation*}
    is a moment-type estimator with asymptotic variance-covariance matrix $\boldsymbol{T}$ such that
    \begin{equation}\label{eq:SA-covariance}
        \boldsymbol{T}_{ij} = \frac{\alpha_i\alpha_j}{(\alpha_0+1)(k-1)^2}s(\ba) + \frac{\alpha_i\alpha_j(\alpha_0+2)-\alpha_0(\alpha_i+\alpha_j)}{(\alpha_0+1)(k-1)} + \frac{\alpha_0\alpha_i}{\alpha_0+1}\delta_{ij}, \quad i,j \in [k],
    \end{equation}
    where $s(\ba):=\sum_{m=1}^k \alpha_m(\alpha_0-\alpha_m)\psi_1(\alpha_m)$ and $\psi_1$ the trigamma function. Its $i$-th diagonal asymptotic variance is
    \begin{equation}\label{eq:SA-diag}
    \tau_i^2(\boldsymbol{\alpha})=\frac{\alpha_i^2}{(\alpha_0+1)(k-1)^2}s(\ba)+\frac{\alpha_i^2(\alpha_0+2)+(k-3)\alpha_0\alpha_i}{(\alpha_0+1)(k-1)}, \quad i \in [k].
    \end{equation}
\end{proposition}

\section{Asymptotic efficiency comparison in the beta distribution case}\label{sec:SA-M1}
\noindent
In this section, we compare the SA and M1 asymptotic variances for $k=2$, i.e., the beta distribution. Throughout the section, we use the common notation $(\alpha, \beta)$ instead of $(\alpha_1, \alpha_2)$. \citet{tamae2020score} proved that  
\[
\tau_i^2(\alpha,\alpha)<\sigma_{i}^2(\alpha,\alpha), \qquad \alpha\in\mathbb{R}_{+},\,\, i=1,2,
\]
that is, they considered the case $\alpha=\beta$. The authors also performed numerical comparisons, in which SA showed lower asymptotic variance than M1, raising the open question of whether SA uniformly dominates M1 over $\mathbb{R}^2_{+}$. We answer this question positively in Theorem \ref{thm:SA-M1-beta}. 

\begin{theorem}\label{thm:SA-M1-beta}
For every $(\alpha,\beta)\in\mathbb{R}^2_{+}$,
\begin{equation*}
\tau_i^2(\alpha,\beta)<\sigma_{i}^2(\alpha,\beta),
\qquad i=1,2.
\end{equation*}
\end{theorem}

\begin{proof}
By symmetry, and without loss of generality, the result is proved for $i = 1$, i.e., for the estimators of $\alpha$. The M1 and SA asymptotic variances, presented in \eqref{eq:M1-diag} and \eqref{eq:SA-diag}, respectively, take the form:
\begin{equation*}
\sigma_{1}^2(\alpha,\beta)
=
\frac{\alpha(\alpha+1)(\alpha+\beta)
\left[2(\alpha+\beta+1)^3-(\alpha+1)\left\{2(\alpha+\beta)^2+\alpha+\beta+1\right\}\right]}
{\beta(\alpha+\beta+1)(\alpha+\beta+2)(\alpha+\beta+3)},
\end{equation*}
and
\begin{equation}\label{eq:SA-beta}
\tau_1^2(\alpha,\beta)
=
\frac{\alpha^3\beta}{\alpha+\beta+1}
\left\{\psi_1(\alpha)+\psi_1(\beta)\right\}
+
\frac{\alpha^2(\alpha+\beta+2)-\alpha(\alpha+\beta)}{\alpha+\beta+1}.
\end{equation}
Let $\overline{\tau}_1^2(\alpha,\beta)$ be obtained from \eqref{eq:SA-beta} by replacing each trigamma term by the upper bound presented in Lemma~\ref{lem:trigamma} of the Appendix: 
\begin{equation}
\overline{\tau}_1^2(\alpha,\beta)
=
\frac{4\alpha^3\beta}{(2\alpha+1)(2\beta+1)}
+
\frac{\alpha(\alpha^2+\beta^2)}{\beta(\alpha+\beta+1)}
+
\frac{\alpha^2(\alpha+\beta+2)-\alpha(\alpha+\beta)}
{\alpha+\beta+1}.
\end{equation}
Then, 
\[
\tau_1^2(\alpha,\beta)<\overline{\tau}_1^2(\alpha,\beta), \quad\forall (\alpha,\beta)\in\mathbb{R}^2_{+},
\]
therefore it suffices to show that 
\[
\overline{\tau}_1^2(\alpha,\beta) <\sigma_{1}^2(\alpha,\beta), \quad\forall (\alpha,\beta)\in\mathbb{R}^2_{+}.
\]
The difference of the two functions takes the form:
\begin{equation}\label{eq:beta-gap}
\sigma_{1}^2(\alpha,\beta)-\overline{\tau}_1^2(\alpha,\beta)
=
\frac{\alpha}
{\beta(\alpha+\beta+1)(\alpha+\beta+2)(\alpha+\beta+3)(2\alpha+1)(2\beta+1)}H(\alpha,\beta),
\end{equation}
where $H$ has the following decomposition:
\begin{equation*}
H(\alpha,\beta)
=5\left(2\alpha^4\beta+\alpha\beta^4-3\alpha^3\beta^2\right)
+2\left(\alpha^4\beta+2\alpha\beta^4-3\alpha^2\beta^3\right)
+10\alpha\beta(\alpha^2-\beta^2)^2
+R(\alpha,\beta),
\end{equation*}
with
\begin{equation*}
\begin{aligned}
R(\alpha,\beta)
={}&4\alpha^5+9\alpha^4\beta+6\alpha^4+2\alpha^3\beta^2+10\alpha^3\beta
+\alpha^2\beta^3+12\alpha^2\beta^2+13\alpha^2\beta+\alpha^2\\
&+20\alpha\beta^4+42\alpha\beta^3
+27\alpha\beta^2+8\alpha\beta+\alpha+4\beta^5+10\beta^4+14\beta^3+7\beta^2+\beta.
\end{aligned}
\end{equation*}
Since the prefactor in \eqref{eq:beta-gap} is positive, it suffices to show that 
\begin{equation}\label{eq:H-positive}
    H(\alpha,\beta)>0, \quad \forall (\alpha,\beta)\in\mathbb{R}^2_{+}.
\end{equation}
It is easy to see that the third term of $H$ is non-negative and every coefficient of $R(\alpha,\beta)$ is positive. The arithmetic--geometric mean inequality gives
\begin{equation*}
2\alpha^4\beta+\alpha\beta^4\geq3\alpha^3\beta^2,
\qquad
\alpha^4\beta+2\alpha\beta^4\geq3\alpha^2\beta^3,
\end{equation*}
therefore the first two terms of $H(\alpha,\beta)$ are also non-negative. Thus, \eqref{eq:H-positive} holds and the result follows.
\end{proof}

\begin{remark}\label{remark:sa-m1-counterexample}
    Theorem~\ref{thm:SA-M1-beta} does not extend uniformly to dimensions $k\geq3$. A simple counterexample for $k=3$ is $\ba=(0.2,5,10)$, which results in
    \begin{equation*}
    \sigma_{2}^2(\boldsymbol{\alpha})\approx49.6613545
    <
    50.1382144\approx\tau_2^2(\boldsymbol{\alpha}).
    \end{equation*}
    The displayed values, rounded to seven decimal places, disprove a uniform componentwise ordering of SA below M1 in $k=3$ dimensions. Similar counterexamples exist for larger $k$.
\end{remark}

\section{A moment estimator based on all components}\label{sec:M2}
\noindent
We now introduce a new moment estimator $\widehat{\boldsymbol{\alpha}}_N$, referred to as M2, which is proportional to the empirical mean vector and equivariant under a common permutation of the coordinates. \citet{ng2011dirichlet} considered a similar estimator, based on pooling moment information from $k-1$ coordinates. With the $k$-th coordinate omitted, their estimator of $\alpha_i$ can be written as
\[
\left\{\frac{\sum_{m=1}^{k-1}\overline X_m(1-\overline X_m)}
{\sum_{m=1}^{k-1}\left(\overline{X_m^2}-\overline{X}_m^{\,2}\right)}
-1\right\}\overline X_i.
\]
This construction depends on the coordinate omitted from the sums. M2 instead pools information across all $k$ coordinates; this removes the need to select an omitted component and ensures permutation equivariance. We formally present M2 in Proposition~\ref{prop:dirichlet-me2} and derive its asymptotic variance-covariance matrix. We adopt the following helpful notation:
\begin{equation*}
\qquad
b(\boldsymbol{\alpha})=\alpha_0^2-\sum_{j=1}^k \alpha_j^2,
\qquad
c(\boldsymbol{\alpha})=\alpha_0^3-\sum_{j=1}^k \alpha_j^3,
\qquad 
\boldsymbol{\alpha}\in\mathbb{R}^k_{+}, \, k\geq2.
\end{equation*}
It is immediate that $b(\boldsymbol{\alpha}), c(\boldsymbol{\alpha})>0, \, \forall\boldsymbol{\alpha}\in\mathbb{R}^k_{+}, \, k\geq2$. 

\begin{lemma}\label{lemma:dirichlet-moments}
If $\mathbf{X}\sim\mathcal{D}_k(\boldsymbol{\alpha})$, $k\geq 2$, then:
\begin{equation*}
\EE\left(\prod_{i=1}^k X_i^{r_i}\right)
=
\frac{\Gamma(\alpha_0)}{\Gamma(\alpha_0+r_0)}
\prod_{i=1}^k\frac{\Gamma(\alpha_i+r_i)}{\Gamma(\alpha_i)},
\qquad
r_0=\sum_{i=1}^k r_i, \quad r_i\geq0.
\end{equation*}
In particular,
\begin{equation*}
\EE(X_i)=\frac{\alpha_i}{\alpha_0},
\qquad
\EE(X_iX_j)
=
\frac{\alpha_i(\alpha_j+\delta_{ij})}{\alpha_0(\alpha_0+1)} \qquad i,j\in [k].
\end{equation*}
\end{lemma}

\begin{proposition} \label{prop:dirichlet-me2}
    Let $\mathbf{X}_n, \, n \in [N], \,N\geq2$, be a random sample from $\mathcal{D}_k(\boldsymbol{\alpha})$. The Moment Estimator
    \begin{equation}\label{eq:M2-estimator}
        \widehat{\alpha}_{N,i}= \frac{1-\sum_{m=1}^k\overline{X_m^2}}{\sum_{m=1}^k\left(\overline{X_m^2}-\overline{X_m}^{\,2}\right)}\overline{X_i}, \qquad i\in[k],
    \end{equation}
    is a moment-type estimator with asymptotic variance-covariance matrix $\boldsymbol{W}$ such that: 
    \begin{equation}\label{eq:M2-covariance}
        \boldsymbol{W}_{ij}
        =\frac{\alpha_0\alpha_i}{\alpha_0+1}\delta_{ij}
        +\frac{2\alpha_0\alpha_i\alpha_j(\alpha_i+\alpha_j)}{b(\boldsymbol{\alpha})(\alpha_0+2)}
        +\frac{\alpha_0\alpha_i\alpha_jp(\ba)}
        {b^2(\boldsymbol{\alpha})(\alpha_0+1)(\alpha_0+2)(\alpha_0+3)}, \quad i,j \in [k],
    \end{equation}
    where
    \begin{equation}\label{eq:p-definition}
    p(\ba)
    =b^2(\boldsymbol{\alpha})(2\alpha_0^2+\alpha_0+1)+
    (\alpha_0-1)\{4c(\boldsymbol{\alpha})(\alpha_0+1)^2-2b(\boldsymbol{\alpha})(\alpha_0+1)(3\alpha_0^2+2\alpha_0+1)\}.
    \end{equation}
    Its $i$-th diagonal asymptotic variance is
    \begin{equation}\label{eq:M2-diag}
        w_{i}^2(\boldsymbol{\alpha}) = \frac{\alpha_0\alpha_i}{\alpha_0+1} + \frac{4\alpha_0\alpha_i^3}{b(\boldsymbol{\alpha})(\alpha_0+2)} + \frac{\alpha_0\alpha_i^2p(\ba)}{b^2(\boldsymbol{\alpha})(\alpha_0+1)(\alpha_0+2)(\alpha_0+3)}, \quad i \in [k].
    \end{equation}
\end{proposition}

\begin{proof} 
    It is immediate to see that Definition \ref{def:me} is satisfied by the functions 
    \[
    h(\mathbf{x}) = (x_1, \dots, x_k, x^2_1, \dots, x^2_k), \quad
    g_i(\mathbf{y})=u(\mathbf{y})y_i, \ i \in [k], \quad u(\mathbf{y}) = \frac{1-\sum_{m=1}^k y_{k+m}} {\sum_{m=1}^k(y_{k+m}-y_m^2)},
    \]
    therefore $\widehat{\boldsymbol{\alpha}}_N$ is indeed a moment-type estimator. Note that $1-\sum_{m=1}^k\overline{X_m^2}\geq0$, with equality iff every observation is a standard basis vector, not necessarily the same one. Similarly, $\sum_{m=1}^k\left(\overline{X_m^2}-\overline{X_m}^{\,2}\right)\geq0$, with equality iff all observations are equal. Since the observations are independent and the Dirichlet distribution is continuous on the open simplex, the two events have probability zero for $N\geq2$. Therefore, $\widehat{\boldsymbol{\alpha}}_N$ takes values in $\mathbb{R}^k_{+}$ almost surely.
    
    The asymptotic variance-covariance matrix of this estimator is of the form $\boldsymbol{W} = \mathbf{G}\mathbf{V}\mathbf{G}^\top$. The $k\times 2k$ Jacobian matrix of $g$ can be broken down into two $k\times k$ matrices:
    \[
        \mathbf{G} = 
        \left[\begin{matrix}
            \mathbf{G}_{1} & \mathbf{G}_{2}
        \end{matrix}\right],
    \]
    with derivatives
    \begin{align*}
       \mathbf{G}_{1;\,ij}(\mathbf{y})&= \frac{\partial g_i}{\partial y_j}(\mathbf{y}) =  \delta_{ij}u(\mathbf{y}) +
        \frac{2y_iy_j\,u(\mathbf{y})}{\sum_{m=1}^{k} (y_{k+m}-y_m^2)},\quad i,j \in [k], \\
        \mathbf{G}_{2;\,ij}(\mathbf{y})&=\frac{\partial g_i}{\partial y_{k+j}}(\mathbf{y})=  \frac{y_i\left(\sum_{m=1}^{k}y_m^2-1\right)}{\left[\sum_{m=1}^{k} (y_{k+m}-y_m^2)\right]^2},
        \quad i,j \in [k].
    \end{align*}
    By evaluating at $\boldsymbol{\mu} = e_h(\boldsymbol{\alpha})$ (see Lemma \ref{lemma:dirichlet-moments}), i.e., for
    \[
    y_i = \frac{\alpha_i}{\alpha_0}, \quad y_{k+i} = \frac{\alpha_i(\alpha_i+1)}{\alpha_0(\alpha_0+1)}, \quad i\in[k],
    \]
    the elements of $\mathbf{G}_1$ and $\mathbf{G}_2$ take the form
    \[
        \mathbf{G}_{1;ij} = \alpha_0\,\delta_{ij} +2b(\boldsymbol{\alpha})^{\,-1}\alpha_0(\alpha_0+1)\,\alpha_i\alpha_j, \quad
        \mathbf{G}_{2;\,ij}=- b^{-1}(\boldsymbol{\alpha})\alpha_0(\alpha_0+1)^2\alpha_i,
        \qquad i,j \in [k].
    \] 
    $\mathbf{G}$ can be represented in matrix form as
    \begin{equation*}
    \mathbf{G}= \left[\, \alpha_0 \mathbf{I}_k + 2 b^{-1}(\boldsymbol{\alpha})\alpha_0(\alpha_0+1)\,\ba \ba^{\top}\qquad - b^{-1}(\boldsymbol{\alpha})\alpha_0(\alpha_0+1)^2 \,\ba\, \boldsymbol{1}_k^{\top}\,\right].
    \end{equation*}
    $\mathbf{V} = v_h(\boldsymbol{\alpha})$ can be broken down to four $k\times k$ blocks:
     \[
        \mathbf{V} = 
        \left[\begin{matrix}
            \mathbf{V}_{11} & \mathbf{V}_{12} \\
            \mathbf{V}_{21} & \mathbf{V}_{22} 
        \end{matrix}\right],
    \] 
    the form of which can be calculated from Lemma \ref{lemma:dirichlet-moments}:
    \begin{eqnarray*}
    \mathbf{V}_{11;ij} &=&  \CC(X_i, X_j) = \frac{(\delta_{ij}\alpha_0-\alpha_j)\,\alpha_i}{\alpha_0^2(\alpha_0+1)},\\
    \mathbf{V}_{12;ij} &=& \CC(X_i, X^2_j)=\frac{2(\delta_{ij}\alpha_0-\alpha_j)\,\alpha_i\,(\alpha_j+1)}{\alpha_0^2(\alpha_0+1)(\alpha_0+2)},\\
    \mathbf{V}_{22;ij} &=& \CC(X^2_i, X^2_j) =\frac{2(\delta_{ij}\alpha_0-\alpha_j)\,\alpha_i\,(\alpha_i+1)}{\alpha_0^2(\alpha_0+1)^2(\alpha_0+2)(\alpha_0+3)}\left\{\delta_{ij}\alpha_0 +(\alpha_j+1)(2\alpha_0+3) \right\},
    \end{eqnarray*}
   for $i,j \in [k]$, and $\mathbf{V}_{21} = \mathbf{V}^\top_{12}$. The equivalent matrix forms are:
    \begin{align*}
    \mathbf{V}_{11} &= c_{11} \left[\alpha_0 \mathbf{D}_{\ba}-\ba\, \ba^{\top}\right], \\
    \mathbf{V}_{12} &= c_{12} \left[\alpha_0 \mathbf{D}_{1}-\ba\, \mathbf{d}_1^{\top}\right], \\
    \mathbf{V}_{22} &= c_{22} \left[-(2\alpha_0+3)\mathbf{d}_1\mathbf{d}_1^{\top}+\alpha_0(\alpha_0+1)\mathbf{D}_2\right],
    \end{align*}
    where 
    \begin{align*}
        \mathbf{D}_{\ba}&=\text{diag}\{\alpha_i\}_{i\in [k]}, \\ 
        \mathbf{D}_1&=\text{diag}\{\alpha_i(\alpha_i+1)\}_{i\in [k]}, \\ 
        \mathbf{D}_2&=\text{diag}\{\alpha_i(\alpha_i+1)(2\alpha_i+3)\}_{i\in [k]}, \\ 
        \mathbf{d}_1&=(\alpha_1(\alpha_1+1), \dots, \alpha_k(\alpha_k+1))^\top,
    \end{align*}
    and
    \begin{equation*}
    c_{11}^{-1} = \alpha_0^2(\alpha_0+1), \quad
    c_{22}^{-1} = \frac{1}{2}\alpha_0^2(\alpha_0+1)^2(\alpha_0+2)(\alpha_0+3), \quad
    c_{12}^{-1} = \frac{1}{2}\alpha_0^2(\alpha_0+1)(\alpha_0+2).
    \end{equation*}
    The block structures of $\mathbf{G}$ and $\mathbf{V}$ simplify the calculations needed to derive the asymptotic variance-covariance matrix: 
    \[
    \boldsymbol{W} = \mathbf{G}\mathbf{V}\mathbf{G}^\top = \sum_{n=1}^2 \sum_{m=1}^2 \mathbf{G}_{n}\mathbf{V}_{nm}\mathbf{G}^\top_{m}.
    \]
    Performing the calculations yields \eqref{eq:M2-covariance}. Taking $i=j$ in \eqref{eq:M2-covariance} gives \eqref{eq:M2-diag}.
\end{proof}

\begin{corollary}\label{cor:M1-M2-beta}
Let $\mathbf{X}_n = (X_{1n}, X_{2n}), \, n \in [N], \, N\geq 2$, be a random sample from $\mathcal{D}_2(\boldsymbol{\alpha})$. Then, $\widetilde{\ba}_N=\widehat{\ba}_N$ almost surely.
\end{corollary}
    
\begin{proof}
    For $k = 2$, $X_{2n} = 1-X_{1n}, \, n \in [N]$, so it is easy to show that:
    \begin{equation*}
        1 - \overline{X_1^2} - \overline{X_2^2} = 2\left(\overline{X_1}-\overline{X_1^2}\right) = 2\left(\overline{X_2}-\overline{X_2^2}\right), \qquad \overline{X_1^2}-\overline{X_1}^2 = \overline{X_2^2}-\overline{X_2}^2.
    \end{equation*}
    For the first shape parameter:
    \begin{equation*}
        \widehat{\alpha}_{N,1}= \frac{1-\sum_{m=1}^k\overline{X_m^2}}{\sum_{m=1}^k\left(\overline{X_m^2}-\overline{X_m}^{\,2}\right)}\overline{X_1} = \frac{2\left(\overline{X_1}-\overline{X_1^2}\right)}{2\left(\overline{X_1^2}-\overline{X_1}^{\,2}\right)}\overline{X_1} = \frac{\overline{X_1}-\overline{X_1^2}}{\overline{X_1^2}-\overline{X_1}^{\,2}}\overline{X_1} = \widetilde{\alpha}_{N,1}.
    \end{equation*}
    The second component follows by symmetry. The above equality holds iff the denominators are non-zero, which is true almost surely.
\end{proof}
As shown in Corollary \ref{cor:M1-M2-beta}, M1 and M2 coincide in the special case of the beta distribution ($k=2$). In Theorem \ref{thm:M2-M1}, we prove that M2 has strictly smaller componentwise asymptotic variances than M1 for $k\geq 3$ dimensions. 
\begin{theorem}\label{thm:M2-M1}
For every $\ba\in\mathbb{R}^k_{+}$, $k\geq3$,
\begin{equation*}
w_{i}^2(\boldsymbol{\alpha})<\sigma_{i}^2(\boldsymbol{\alpha}), \quad i\in[k].
\end{equation*}
\end{theorem}

\begin{proof}
    The difference of the two variances can be rewritten as $\sigma_{i}^2(\boldsymbol{\alpha})-w_{i}^2(\boldsymbol{\alpha}) = \pi_1(\boldsymbol{\alpha})\pi_2(\boldsymbol{\alpha})$, where
    \begin{align*}
        \pi_1(\boldsymbol{\alpha}) &= \frac{\alpha_0\alpha_i}{b^2(\boldsymbol{\alpha})(\alpha_0+2)(\alpha_0+3)(\alpha_0-\alpha_i)}, \\[2ex]
        \pi_2(\boldsymbol{\alpha}) &= b^2(\boldsymbol{\alpha}) \left[2(\alpha_0+3)\alpha_i+(1-\alpha_0)^2\right] + 4\alpha_i(\alpha_0-\alpha_i)(1-\alpha_0)(1+\alpha_0)c(\boldsymbol{\alpha}) \\
        &\quad -2b(\boldsymbol{\alpha})(\alpha_0-\alpha_i)\alpha_i\left[2(\alpha_0+3)\alpha_i +(1-\alpha_0)^2+3\alpha_0(1-\alpha_0)(1+\alpha_0)\right].
    \end{align*}
    It is immediate that 
    \[
    \pi_1(\boldsymbol{\alpha}) > 0, \quad \forall \boldsymbol{\alpha}\in\mathbb{R}^k_{+}, \, k\geq 3,
    \]
    so it suffices to show the same for $\pi_2(\boldsymbol{\alpha})$. After calculations, $\pi_2(\boldsymbol{\alpha})$ can be expressed as a 4th degree polynomial of $\alpha_i$:
    \[
    \pi_2(\boldsymbol{\alpha}) = r_4(\boldsymbol{\alpha}_{-i})\alpha^4_i+r_3(\boldsymbol{\alpha}_{-i})\alpha^3_i+r_2(\boldsymbol{\alpha}_{-i})\alpha^2_i+r_1(\boldsymbol{\alpha}_{-i})\alpha_i+r_0(\boldsymbol{\alpha}_{-i}),
    \]
    with $\boldsymbol{\alpha}_{-i}:=(\alpha_1, \dots, \alpha_{i-1} , \alpha_{i+1}, \dots, \alpha_k)$ and
    \begin{align*}
        r_4(\boldsymbol{\alpha}) &= 6 \alpha_0b(\boldsymbol{\alpha}), \\[1.5ex]
        r_3(\boldsymbol{\alpha}) &= 2\alpha_0\left\{d(\boldsymbol{\alpha})+3(2\alpha_0+1)b(\boldsymbol{\alpha})\right\}, \\[1.5ex]
        r_2(\boldsymbol{\alpha}) &= (6\alpha^3_0+8\alpha^2_0+2\alpha_0)b(\boldsymbol{\alpha})+3b^2(\boldsymbol{\alpha})+4\alpha^2_0d(\boldsymbol{\alpha}), \\[1.5ex]
        r_1(\boldsymbol{\alpha}) &= 4(\alpha_0+1)b^2(\boldsymbol{\alpha})+2\alpha_0b(\boldsymbol{\alpha})(\alpha_0-1)^2+2\alpha_0(\alpha^2_0-1)d(\boldsymbol{\alpha}), \\[1.5ex]
        r_0(\boldsymbol{\alpha}) &= (1-\alpha_0)^2b^2(\boldsymbol{\alpha}),
    \end{align*}
    where $d(\boldsymbol{\alpha}):=3\alpha_0b(\boldsymbol{\alpha})-2c(\boldsymbol{\alpha})$. It is immediate to see that $r_{4}(\boldsymbol{\alpha}), r_{3}(\boldsymbol{\alpha}), r_{2}(\boldsymbol{\alpha})>0$ and $r_{0}(\boldsymbol{\alpha})\geq0$ for $k\geq 2$, therefore the same holds for these coefficients evaluated at $\boldsymbol{\alpha}_{-i}$ for $k\geq 3$. The only case that is not immediate is $r_1(\boldsymbol{\alpha})$ when $\alpha_0<1$. By Lemma \ref{lem:beta-gamma-bound} of the Appendix,
    \begin{align*}
    r_1(\boldsymbol{\alpha})&>
    4(\alpha_0+1)b^2(\boldsymbol{\alpha})+2\alpha_0b(\boldsymbol{\alpha})(\alpha_0-1)^2+2b^2(\boldsymbol{\alpha})(\alpha_0^2-1) \\
    &=2(\alpha_0+1)^2b^2(\boldsymbol{\alpha})+2\alpha_0b(\boldsymbol{\alpha})(\alpha_0-1)^2>0.
    \end{align*}
    Therefore, $\pi_2(\boldsymbol{\alpha}) > 0, \, \forall \boldsymbol{\alpha}\in\mathbb{R}^k_{+}, \, k\geq 3$, and the result follows.
\end{proof}
\noindent
Corollary~\ref{cor:M2-M1} shows that both moment estimators have positive asymptotic covariances.
\begin{corollary}\label{cor:M2-M1}
For every $\ba\in\mathbb{R}^k_{+}$, $k\geq2$,
\begin{equation*}
\boldsymbol{\Sigma}_{ij}(\boldsymbol{\alpha})>0, \quad \boldsymbol{W}_{ij}(\boldsymbol{\alpha})>0, \quad i,j\in[k].
\end{equation*}
\end{corollary}
\begin{proof}
    The proof focuses on the covariances ($i\neq j$), as the variance case ($i=j$) is immediate. For M1, the result follows from the definition of $\boldsymbol{\Sigma}_{ij}$ in \eqref{eq:Sij}:
    \begin{equation*}
    \boldsymbol{\Sigma}_{ij} = \frac{\alpha_i\alpha_j\alpha_0(\alpha_i + 1)(\alpha_j + 1)(2\alpha_0^2+\alpha_0+1)}{(\alpha_0-\alpha_i)(\alpha_0-\alpha_j)(\alpha_0+1)(\alpha_0+2)(\alpha_0+3)}>0, \quad i,j \in [k], i\neq j.
    \end{equation*}
    For M2, the result follows from the definition of $\boldsymbol{W}_{ij}$ in \eqref{eq:M2-covariance} and Lemma \ref{lemma:p-positive} of the Appendix that establishes the positivity of $p(\ba)$:
    \begin{equation*}
        \boldsymbol{W}_{ij}
        =\frac{2\alpha_0\alpha_i\alpha_j(\alpha_i+\alpha_j)}{b(\boldsymbol{\alpha})(\alpha_0+2)}
        +\frac{\alpha_0\alpha_i\alpha_jp(\ba)}
        {b^2(\boldsymbol{\alpha})(\alpha_0+1)(\alpha_0+2)(\alpha_0+3)}>0, \quad i,j \in [k], i\neq j.
    \end{equation*}
\end{proof}

\section{Asymptotic efficiency comparison in the Dirichlet case}\label{sec:SA-M2}

\noindent
We now compare SA and M2 in the general Dirichlet case (dimension $k\geq2$). Lemma \ref{lem:SA-M2-equal} shows that SA outperforms M2 in the symmetric case where all the shape parameters are equal. This extends the corresponding result of \citet{tamae2020score} for the specific case of $k=2$, i.e., the beta distribution. Furthermore, in Theorem~\ref{thm:SA-M2-matrix-comparison} we establish the stronger result that SA is asymptotically more efficient than M2 when all shape parameters are equal and that, in any other case, the difference $\mathbf W(\boldsymbol\alpha)-\mathbf T(\boldsymbol\alpha)$ is indefinite, so neither estimator dominates the other in terms of full asymptotic efficiency.

\begin{lemma}\label{lem:SA-M2-equal}
Let $k\geq2$ and $\ba=a\mathbf{1}_k$, $a>0$. Then
\begin{equation*}
\tau_i^2(\boldsymbol{\alpha})<w_{i}^2(\boldsymbol{\alpha}), \quad i\in[k].
\end{equation*}
\end{lemma}

\begin{proof}
For $\ba=a\mathbf1_k$,
\[
\alpha_0=ka,\qquad b(\ba)=k(k-1)a^2,\qquad c(\ba)=k(k-1)(k+1)a^3.
\]
Substituting into the asymptotic variance-covariance matrices of SA and M2, presented in \eqref{eq:SA-covariance} and \eqref{eq:M2-covariance}, respectively, and collecting terms gives
\begin{equation}
\begin{aligned}
   \mathbf T(\ba)&=\frac{ka^2}{ka+1}\mathbf I_k+\frac{a^2\{ka^2\psi_1(a)+ka-2k+2\}}{(k-1)(ka+1)}\mathbf1_k\mathbf1_k^\top, \\
   \mathbf W(\ba)&=\frac{ka^2}{ka+1}\mathbf I_k+\left\{\frac{2a(a+1)(ka+1)^2}{(k-1)(ka+2)(ka+3)}-\frac{a^2}{ka+1}\right\}\mathbf1_k\mathbf1_k^\top.
\end{aligned}    
\end{equation}
The identity-matrix terms cancel, yielding
\[
\mathbf D(\ba):=\mathbf W(\ba)-\mathbf T(\ba)=\delta\mathbf1_k\mathbf1_k^\top,
\]
where
\[
\delta=\frac{2a(a+1)(ka+1)^2}{(k-1)(ka+2)(ka+3)}-\frac{a^2\{ka^2\psi_1(a)+ka-k+1\}}{(k-1)(ka+1)}.
\]
In particular, $w_i^2(\ba)-\tau_i^2(\ba)=\delta$ for every $i\in[k]$. By Lemma~\ref{lem:trigamma} of the Appendix,
\[
\psi_1(a)<\frac{1}{a+\tfrac12}+\frac{1}{a^2}.
\]
Since the coefficient of $\psi_1(a)$ in $\delta$ is negative, substitution of this bound and simplification give
\[
\delta>\frac{a\{a(k-2)+1\}\left(5k^2a^3+2k^2a^2+11ka^2+4ka+4a+2\right)}{(2a+1)(k-1)(ka+1)(ka+2)(ka+3)}>0.
\]
The final inequality follows from $a>0$ and $k\ge2$, proving the result.
\end{proof}

\begin{theorem}\label{thm:SA-M2-matrix-comparison}
Let $k\ge2$ and
\[
\mathbf D(\ba):=\mathbf W(\boldsymbol\alpha)-\mathbf T(\boldsymbol\alpha), \qquad \ba\in\mathbb{R}^k_+.
\]
If there exists $a\in\mathbb{R}_+$ such that $\ba=a\mathbf{1}_k$, then $\mathbf D(\ba)\succeq0$ and has rank one. Otherwise, $\mathbf D(\ba)$ has exactly one positive eigenvalue, one negative eigenvalue, and $k-2$ zero eigenvalues.
\end{theorem}

\begin{proof}
Subtracting the covariance formulas for SA and M2, presented in \eqref{eq:SA-covariance} and \eqref{eq:M2-covariance} respectively, gives
\begin{equation}\label{eq:Da}
    \mathbf D(\ba)=\mathbf W(\ba)-\mathbf T(\ba)=\ba\mathbf d^\top+\mathbf d\ba^\top+\gamma\ba\ba^\top,
\end{equation}
where $\gamma\in\mathbb{R}$ and
\[
A=\frac{2\alpha_0}{b(\ba)(\alpha_0+2)},\qquad B=\frac{\alpha_0}{(\alpha_0+1)(k-1)}, \qquad d_i=A\alpha_i^2+B,\,\, i\in[k].
\]
The columns of $\mathbf D$ belong to $\operatorname{span}\{\ba,\mathbf d\}$, so $\operatorname{rank}(\mathbf D)\le2$. 

Suppose that the shape parameters are all equal, i.e., there exists $a>0$ such that $\ba=a\mathbf1_k$.  Then, $\mathbf D=\delta\mathbf1_k\mathbf1_k^\top$
by Lemma~\ref{lem:SA-M2-equal}. Consequently, for every $\mathbf z\in\mathbb R^k$,
\[
\mathbf z^\top\mathbf D\mathbf z=\delta(\mathbf1_k^\top\mathbf z)^2\ge0.
\]
Therefore, $\mathbf D$ is positive semidefinite and has rank one.

Suppose that the shape parameters are not all equal, and choose $i,j$ such that
\[
\alpha_i=\min_{r\in[k]}\alpha_r < \max_{r\in[k]}\alpha_r=\alpha_j.
\]
Since
\[
b(\ba)=2\sum_{r<t}\alpha_r\alpha_t\ge2\alpha_j\sum_{r\ne j}\alpha_r\ge2(k-1)\alpha_i\alpha_j,
\]
it follows that
\[
\frac BA=\frac{b(\ba)(\alpha_0+2)}{2(k-1)(\alpha_0+1)}\ge \alpha_i\alpha_j\frac{\alpha_0+2}{\alpha_0+1}>\alpha_i\alpha_j.
\]
The symmetric principal submatrix
\[
\mathbf H=\begin{pmatrix}\mathbf D_{ii}& \mathbf D_{ij}\\\mathbf D_{ij}&\mathbf D_{jj}\end{pmatrix},
\]
has
\[
\det(\mathbf H)= D_{ii}D_{jj}-D_{ij}^2=-(\alpha_i d_j-\alpha_j d_i)^2=-(\alpha_j-\alpha_i)^2(A\alpha_i\alpha_j-B)^2<0.
\]
Consequently, $\mathbf H$ has two eigenvalues of opposite signs, since their product is $\det(\mathbf H)<0$. Hence, there exist $\mathbf u_+,\mathbf u_-\in\mathbb R^2$ such that
\[
\mathbf u_+^\top\mathbf H\mathbf u_+>0,\qquad \mathbf u_-^\top\mathbf H\mathbf u_-<0.
\]
Embedding these vectors into $\mathbb R^k$ by placing their entries in positions $i,j$ and zeros elsewhere gives vectors $\mathbf z_+,\mathbf z_-$ satisfying
\[
\mathbf z_+^\top\mathbf D\mathbf z_+=\mathbf u_+^\top\mathbf H\mathbf u_+>0,\qquad \mathbf z_-^\top\mathbf D\mathbf z_-=\mathbf u_-^\top\mathbf H\mathbf u_-<0.
\]
Thus, the quadratic form of $\mathbf D$ takes both positive and negative values. Since $\mathbf D$ is symmetric, it has at least one positive and at least one negative eigenvalue. Together with $\operatorname{rank}(\mathbf D)\le2$, this shows that it has exactly one eigenvalue of each sign and $k-2$ zero eigenvalues, concluding the proof.
\end{proof}

\noindent
Theorem \ref{thm:SA-M2-matrix-comparison} shows that the difference between the asymptotic covariance matrices of M2 and SA is indefinite for unequal parameters. We now give a sufficient condition for all its diagonal entries to be positive.

\begin{theorem}\label{thm:SA-M2-c1}
Let $k\geq3$. If $\ba\in[1, +\infty)^k$, then
\begin{equation*}
\tau_i^2(\ba)<w_{i}^2(\ba), \qquad i\in[k].
\end{equation*}
\end{theorem}

\begin{proof}
Let $\ba\in[1,+\infty)^k$. From \eqref{eq:SA-diag} and \eqref{eq:M2-diag}, the difference of the two estimators' asymptotic variances can be written as
\begin{equation*}
w_i^2(\ba)-\tau_i^2(\ba)=\frac{\alpha_i}{(k-1)^2b^2(\ba)(\alpha_0+1)(\alpha_0+2)(\alpha_0+3)}M_i(\ba),
\end{equation*}
where the prefactor is positive and
\begin{equation*}
M_i(\ba)=U(\ba)\alpha_i^2+V(\ba)\alpha_i+Z(\ba),
\end{equation*}
with
\begin{align*}
U(\ba)&=4(k-1)^2b(\ba)\alpha_0(\alpha_0+1)(\alpha_0+3),\\
V(\ba)&=(k-1)^2\alpha_0p(\ba)-b^2(\ba)(\alpha_0+2)(\alpha_0+3)\left[s(\ba)+(k-1)(\alpha_0+2)\right],\\
Z(\ba)&=2(k-1)\alpha_0b^2(\ba)(\alpha_0+2)(\alpha_0+3).
\end{align*}
Since $U(\ba)>0$ and $Z(\ba)>0$, completing the square gives the identity
\begin{equation}\label{eq:Mi}
M_i(\ba) = \left\{\sqrt{U(\ba)}\,\alpha_i-\sqrt{Z(\ba)}\right\}^2+\alpha_i\left\{V(\ba)+2\sqrt{U(\ba)Z(\ba)}\right\}.
\end{equation}
The squared term in \eqref{eq:Mi} is non-negative. Hence, proving that
\begin{equation*}
Q(\ba):=V(\ba)+2\sqrt{U(\ba)Z(\ba)}>0, \qquad \ba\in[1,+\infty)^k,
\end{equation*}
is sufficient to prove the result. Expanding $Q(\ba)$ results in
\begin{equation}\label{eq:Q1}
\begin{aligned}
Q(\ba)
={}&
(k-1)^2\alpha_0b^2(\ba)(2\alpha_0^2+\alpha_0+1)-(k-1)b^2(\ba)(\alpha_0+2)^2(\alpha_0+3) \\
&-2(k-1)^2\alpha_0b(\ba)(\alpha_0-1)(\alpha_0+1)(3\alpha_0^2+2\alpha_0+1)\\
&+4(k-1)^2\alpha_0c(\ba)(\alpha_0+1)^2(\alpha_0-1)\\
&+4(k-1)\alpha_0b(\ba)(\alpha_0+3)\sqrt{2(k-1)b(\ba)(\alpha_0+1)(\alpha_0+2)} \\
&-b^2(\ba)(\alpha_0+2)(\alpha_0+3)s(\ba).
\end{aligned}
\end{equation}
By Lemma~\ref{trigamma-inequality-lemma} of the Appendix,
\begin{equation}\label{eq:Q-trigamma-sum}
s(\ba)=\sum_{m=1}^k\alpha_m(\alpha_0-\alpha_m)\psi_1(\alpha_m)
<\sum_{m=1}^k(\alpha_0-\alpha_m)\left(1+\frac{2}{3\alpha_m}\right)
=(k-1)\alpha_0+\frac23(\alpha_0I-k),
\end{equation}
where $I:=\sum_{i=1}^k\alpha_i^{-1}$. Furthermore,
\begin{equation}\label{eq:b-bound}
    (k-1)(2\alpha_0-k)\le b(\ba)\le\frac{(k-1)\alpha_0^2}{k},
\end{equation}
where the lower bound follows from $\ba\in[1, +\infty)^k$ and the upper bound follows from the Cauchy–Schwarz inequality. Utilizing \eqref{eq:b-bound} and $\alpha_0\ge k$,
\begin{equation}\label{eq:Q-root}
\alpha_0\sqrt{2(k-1)b(\ba)(\alpha_0+1)(\alpha_0+2)}
\ge b(\ba)\sqrt{2k(\alpha_0+1)(\alpha_0+2)}
>\frac76 b(\ba)(2\alpha_0+3).
\end{equation}
Substituting \eqref{eq:Q-trigamma-sum} and \eqref{eq:Q-root} into \eqref{eq:Q1} results in
\begin{equation}\label{eq:Q-bound1}
\begin{aligned}
Q(\ba)
>{}&
(k-1)^2\alpha_0b^2(\ba)(2\alpha_0^2+\alpha_0+1)-(k-1)b^2(\ba)(\alpha_0+2)^2(\alpha_0+3) \\
&-2(k-1)^2\alpha_0b(\ba)(\alpha_0-1)(\alpha_0+1)(3\alpha_0^2+2\alpha_0+1)\\
&+4(k-1)^2\alpha_0c(\ba)(\alpha_0+1)^2(\alpha_0-1)\\
&+\frac{14}{3}(k-1)b^2(\ba)(\alpha_0+3)(2\alpha_0+3)-(k-1)\alpha_0b^2(\ba)(\alpha_0+2)(\alpha_0+3)\\
&-\frac23\alpha_0I b^2(\ba)(\alpha_0+2)(\alpha_0+3)
+\frac{2}{3}kb^2(\ba)(\alpha_0+2)(\alpha_0+3).
\end{aligned}
\end{equation}
Replacing the terms involving $c(\ba)$ and $I$ in \eqref{eq:Q-bound1} by their lower bound in Lemma~\ref{lem:Q-mixed} of the Appendix, and collecting terms results in
\begin{equation}\label{eq:Q-bound2}
\begin{aligned}
Q(\ba)
>{}&(k-1)b^2(\ba)
\left[\frac{\alpha_0}{3}\left\{(k+1)\alpha_0+45-7k\right\}-4k+\frac{2k}{\alpha_0}\right]\\
&+d_k(\alpha_0)
b(\ba)\left[(k-1)\alpha_0^2-kb(\ba)\right],
\end{aligned}
\end{equation}
where
\[
d_k(\alpha_0):=\frac{2(k-1)(\alpha_0-1)^2(\alpha_0+1)}{\alpha_0}
-\frac{2k}{3}(\alpha_0+2)(\alpha_0+3).
\]
To simplify \eqref{eq:Q-bound2}, first note that $3\le k\le\alpha_0$ implies
\[
(k+1)\alpha_0+9-7k
\ge (k+1)k+9-7k
=(k-3)^2\ge0.
\]
Consequently,
\begin{equation}\label{eq:Q2}
\frac{\alpha_0}{3}\left\{(k+1)\alpha_0+45-7k\right\}
\ge12\alpha_0.    
\end{equation}
Moreover, since $2k/[3(k-1)]\le1$,
\begin{equation}\label{eq:Q3}
    \begin{aligned}
        d_k(\alpha_0)&\ge
        (k-1)\left\{
        \frac{2(\alpha_0-1)^2(\alpha_0+1)}{\alpha_0}
        -(\alpha_0+2)(\alpha_0+3)
        \right\}\\
        &=
        (k-1)\left\{
        \left(\alpha_0-\frac72\right)^2-\frac{81}{4}
        +\frac{2}{\alpha_0}
        \right\}
        \ge-24(k-1).
    \end{aligned}
\end{equation}
The upper bound in \eqref{eq:b-bound} gives
$(k-1)\alpha_0^2-kb(\ba)\ge0$; substituting this along with \eqref{eq:Q2} and \eqref{eq:Q3} into \eqref{eq:Q-bound2} gives
\begin{equation}\label{eq:Q-bound3}
    Q(\ba)\ge{}(k-1)b^2(\ba)\left(12\alpha_0-4k+\frac{2k}{\alpha_0}\right)-24(k-1)b(\ba)\left[(k-1)\alpha_0^2-kb(\ba)\right].
\end{equation}
For the remaining term, the lower bound in \eqref{eq:b-bound} gives
\[
2\left[(k-1)\alpha_0^2-kb(\ba)\right]
\le b(\ba)(\alpha_0-k),
\]
therefore \eqref{eq:Q-bound3} simplifies into
\begin{align*}
Q(\ba)
\ge{}&(k-1)b^2(\ba)
\left[12\alpha_0-4k+\frac{2k}{\alpha_0}-12(\alpha_0-k)\right]=(k-1)b^2(\ba)\left(8k+\frac{2k}{\alpha_0}\right)>0.
\end{align*}
Thus, $Q(\ba)>0$ for every $\ba\in[1,+\infty)^k$, and the result follows.
\end{proof}

\begin{remark}\label{remark:SA-M2-counterexample}
    SA does not dominate M2 uniformly over the parameter set. The reverse ordering of the asymptotic variances can also occur. The counterexample in Remark \ref{remark:sa-m1-counterexample} for the SA--M1 comparison also covers the SA--M2 comparison, since M2 has uniformly smaller componentwise asymptotic variances than M1. Hence, neither SA nor M2 has uniformly smaller componentwise asymptotic variances than the other one on $\mathbb{R}^k_{+}, \, k\geq3$. 

    Furthermore, even if $\ba\in[1, +\infty)^k$, the result of Theorem~\ref{thm:SA-M2-c1} does not extend to linear combinations of the parameters. Consider for example the shape vector $\ba=(1,2,2)$. By Theorem~\ref{thm:SA-M2-c1}, SA has smaller asymptotic variance than M2 for each of the three shape parameters. However, if we are interested in estimating the difference $\alpha_1-\alpha_2$, the ordering is reversed. Equations \eqref{eq:SA-covariance} and \eqref{eq:M2-covariance} give asymptotic variances of $31/12+\pi^2/9$ for SA and $175/48$ for M2, approximately $3.680$ and $3.646$, respectively.

\end{remark}

\section{Bounds on the uniform threshold}\label{sec:numerical}
\noindent
For each dimension $k\geq3$, define the optimal uniform threshold
\begin{equation}\label{eq:uniform-threshold}
 c_k^\ast := \inf\left\{c>0:
 w_{i}^2(\ba)-\tau_i^2(\ba)>0
 \text{ for every }i\in[k]
 \text{ and every }\ba\in[c,\infty)^k
 \right\}.
\end{equation}
In this section, we derive mathematical and computational lower bounds for $c_k^\ast$. By Theorem~\ref{thm:SA-M2-c1}, SA has lower asymptotic componentwise variances than M2 when all parameters satisfy $\alpha_i\geq u=1, \, i\in [k]$; therefore, $c_k^\ast\leq u$. 

\subsection{Analytical Bounds}

\begin{lemma}\label{lemma:large-concentration-limit}
Let $k\geq3$ and $c>0$. Consider $\ba_L\in\mathbb{R}_{+}^k$ such that
\begin{equation}\label{eq:large-L-vector-supp}
\ba_L=(
\underbrace{L,\ldots,L}_{m\ \mathrm{times}},
\underbrace{c,\ldots,c}_{k-m\ \mathrm{times}}),\qquad L>c, \, m\in\{2, \dots, k-1\}.
\end{equation}
Then,
\begin{equation}\label{eq:large-concentration-limit}
\lim_{L\to\infty}
\frac{w_{1}^2(\ba_L)-\tau_1^2(\ba_L)}{L^2}
=\frac{2}{m-1} - \frac{(k-m)c\psi_1(c)+m+k-2}{(k-1)^2}.
\end{equation}
\end{lemma}

\begin{proof}
The specific form of $\ba_L$ simplifies the following quantities:
\begin{equation}
\begin{aligned}
\alpha_0&=mL+o(L), \\
b(\boldsymbol{\alpha}_L)&=m(m-1)L^2 + o\left(L^2\right), \\
c(\boldsymbol{\alpha}_L)&=m(m^2-1)L^3 + o\left(L^3\right), \\
p(\boldsymbol{\alpha}_L)&=\left\{2m^4(m-1)^2-6m^5(m-1)+4m^4(m^2-1)\right\}L^6 + o\left(L^6\right).
\end{aligned}
\label{eq:simplified-notation-aL}
\end{equation}
Direct substitution of \eqref{eq:simplified-notation-aL} in $w_{1}^2(\ba)$, presented in \eqref{eq:M2-diag}, results in:
\begin{align*}
    w_{1}^2(\ba_L) &= \frac{\alpha_0\alpha_1}{\alpha_0+1} + \frac{4\alpha_0\alpha_1^3}{b(\boldsymbol{\alpha}_L)(\alpha_0+2)} + \frac{\alpha_0\alpha_1^2p(\ba_L)}{b^2(\boldsymbol{\alpha}_L)(\alpha_0+1)(\alpha_0+2)(\alpha_0+3)} \\
    &= \left\{L+o(L)\right\} + \left\{\frac{4}{m(m-1)}L+o(L)\right\} + \left\{\frac{2}{m-1}L^2+o\left(L^2\right)\right\} \\
    &=\frac{2}{m-1}L^2+o\left(L^2\right).
\end{align*}
Moving on to $\tau^2_1$, Lemma \ref{lem:trigamma} of the Appendix gives
\begin{equation*}
\frac1L<\psi_1(L)
<\frac1L+\frac1{L^2},
\end{equation*}
therefore $L\psi_1(L)=1 + o(1)$. The trigamma sum in $\tau^2_1(\ba)$ takes the form
\begin{equation}
\begin{aligned}
\sum_{r=1}^k\alpha_r(\alpha_0-\alpha_r)\psi_1(\alpha_r)
={}&(k-m)\left\{mL+o(L)\right\}c\psi_1(c) + m  \left\{(m-1)L+o(L)\right\}L\psi_1(L)\\
={}&\left\{m-1+(k-m)c\psi_1(c)\right\}mL+o(L).
\end{aligned}
\label{eq:simplified-SA-sum}
\end{equation}
Direct substitution of \eqref{eq:simplified-notation-aL} and \eqref{eq:simplified-SA-sum} in $\tau_{1}^2(\ba)$, presented in \eqref{eq:SA-diag}, results in
\begin{align*}
    \tau_1^2(\boldsymbol{\alpha}_L)
    &=
    \frac{\alpha_1^2}{(\alpha_0+1)(k-1)^2}
    \sum_{m=1}^k\alpha_m(\alpha_0-\alpha_m)\psi_1(\alpha_m)
    +\frac{\alpha_1^2(\alpha_0+2)+(k-3)\alpha_0\alpha_1}
    {(\alpha_0+1)(k-1)} \\
    &=\left\{\frac{m-1+(k-m)c\psi_1(c)}{(k-1)^2}L^2+o\left(L^2\right)\right\} + \left\{\frac{1}{k-1}L^2+o\left(L^2\right)\right\} \\
    &= \frac{(k-m)c\psi_1(c)+m+k-2}{(k-1)^2}L^2+o\left(L^2\right).
\end{align*}
Thus,
\begin{equation*}
    w_{1}^2(\ba_L)-\tau_1^2(\ba_L)
=\left\{\frac{2}{m-1} - \frac{(k-m)c\psi_1(c)+m+k-2}{(k-1)^2}\right\}L^2 + o\left(L^2\right),
\end{equation*}
and the result follows by taking the limit $L\rightarrow+\infty$.
\end{proof}

\begin{corollary} \label{cor:ck-bounds}
    Let $k\geq 3$. Then, 
    \begin{equation}\label{eq:threshold-bracket}
     \ell_k\leq c_k^\ast\leq1,
    \end{equation}
    where $\ell_k$ is the unique positive solution of
    \begin{equation}\label{eq:lk-definition}
     x\psi_1(x)=3+\frac{2}{k-2}.
    \end{equation}    
\end{corollary}
\begin{proof}
    By Lemma~\ref{lemma:large-concentration-limit},
    \[
    \lim_{L\to\infty}
    \frac{w_{1}^2(\ba_L)-\tau_1^2(\ba_L)}{L^2}
    =\frac{2}{m-1} - \frac{(k-m)c\psi_1(c)+m+k-2}{(k-1)^2}.
    \]
    Lemma~\ref{lem:x-trigamma-decreasing} proves that the function $x\mapsto x\psi_1(x)$ is strictly decreasing from $+\infty$ to $1$. Hence, for fixed $k$ and $m$, the equation
    \begin{equation}\label{eq:boundary-lk}
        \frac{2}{m-1} - \frac{(k-m)x\psi_1(x)+m+k-2}{(k-1)^2} = 0 \quad\Longleftrightarrow\quad x\psi_1(x)=1+\frac{2(k-1)}{m-1},
    \end{equation}
    has a unique positive solution $x=\ell_k(m)$. By the definition of $c^\ast_k$ and Theorem \ref{thm:SA-M2-c1}, it follows that 
    \begin{equation*}
     \ell_k(m)\leq c_k^\ast\leq u=1.
    \end{equation*}
    It is easy to show that $\ell_k(m)$ is strictly increasing in $m\in\{2, \ldots, k-1\}$; therefore, the tightest lower bound is achieved at $m=k-1$. Substitution in \eqref{eq:boundary-lk} results in the simplified equation
    \begin{equation*}
     x\psi_1(x)=3+\frac{2}{k-2},
    \end{equation*}
    concluding the proof.
\end{proof}

\begin{remark}
    Corollary~\ref{cor:ck-bounds} shows that the tightest lower bound within the family of vectors $\ba_L$ considered in Lemma~\ref{lemma:large-concentration-limit} is achieved by vectors of the form 
    \begin{equation*}
        \ba_L=(
        \underbrace{L,\ldots,L}_{k-1\ \mathrm{times}},c),\quad L>c,
    \end{equation*}
    i.e., when one parameter remains fixed, while all the others diverge. The limiting value $\ell_\infty$ is the positive solution of 
    \[
    \ell_\infty\psi_1(\ell_\infty)=3,
    \]
    which yields $\ell_\infty\approx 0.384665$. 
\end{remark}
\subsection{Numerical Bounds}
\noindent
We seek to improve the lower bound $\ell_k$ by computationally finding parameters $\boldsymbol{\alpha}\in[c,\infty)^k$, with $c\geq\ell_k$, for which
\[
\Pi(\boldsymbol{\alpha})=w_1^2(\ba)-\tau_1^2(\ba)<0.
\]
Specifically, let $\mathcal{A}_k$ be the set of all vectors evaluated in the numerical experiment. We define the numerical lower bound $s_k\in[\ell_k,c^\ast_k]$ as the largest candidate $c$ for which a counterexample is found, i.e.,
\[
s_k := \max\left(\{\ell_k\}\cup\left\{\min_{1\le j\le k}\alpha_j :\boldsymbol{\alpha}\in\mathcal{A}_k,\ \Pi(\boldsymbol{\alpha})<0\right\}\right).
\]

Our procedure consisted of two nested searches. The outer search considered $100$ candidate thresholds $c\in[\ell_k+10^{-7}, 0.5]$, with logarithmically spaced distances from $\ell_k$ to concentrate the search near the analytical bound. Each candidate was investigated independently, so an unsuccessful search did not exclude larger thresholds from consideration. For each candidate $c$, the inner search performed $20$ independent generalized simulated-annealing runs over $\ba\in[c,10^6]^k$. All coordinates were free to vary, and the starting coordinates were generated independently. 

Each run was subject to two computational caps: $10{,}000$ outer annealing cycles and approximately $50{,}000$ objective-function evaluations. An outer cycle includes multiple proposal steps, each of which evaluates a candidate parameter vector and then accepts or rejects it; the local-search procedure may also require additional evaluations. Consequently, the number of objective evaluations generally exceeds the number of outer cycles. A run terminated when either the computational limit was reached or, earlier, when a parameter vector satisfying the counterexample criterion was found.

\begin{figure}[t]
    \centering
    \includegraphics[width=\textwidth]{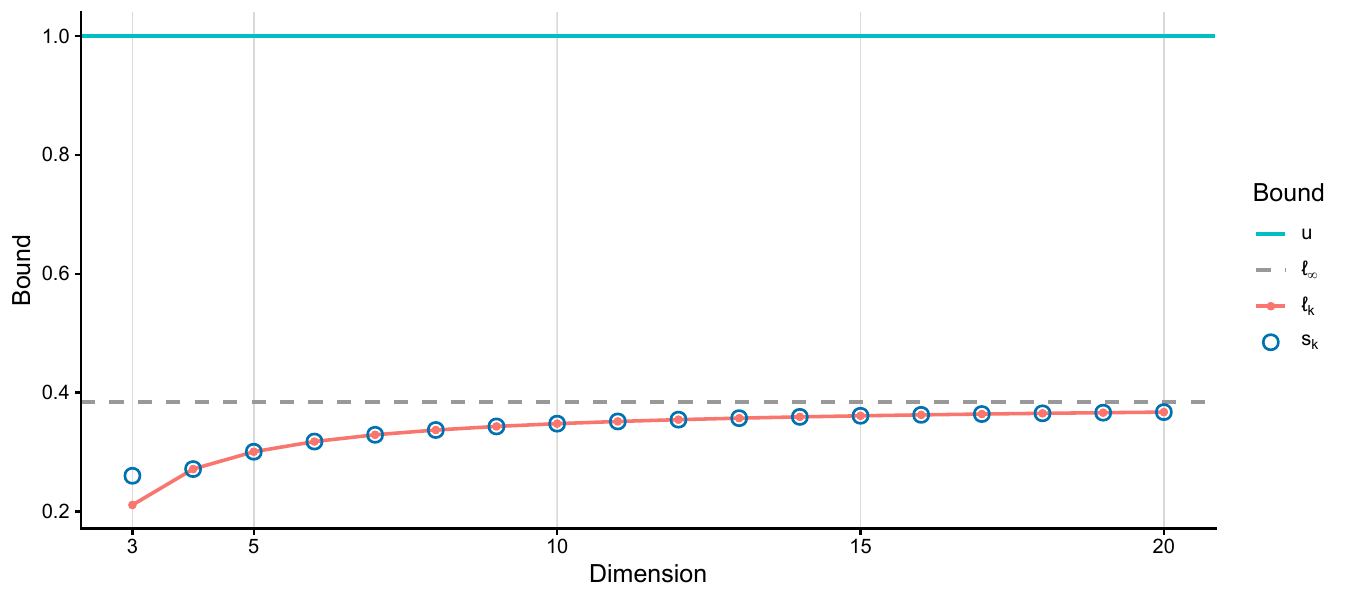}
    \caption{Bounds on the optimal uniform threshold $c_k^\ast$ for dimensions $k\in\{3,\ldots,20\}$. The analytical lower bound $\ell_k$ (red line), defined as the root of~\eqref{eq:lk-definition}, increases toward $\ell_\infty\approx0.384665$ (gray dashed line). The numerical lower bound $s_k$ (blue circles) is derived by looking for numerical counterexamples that improve upon $\ell_k$. An improvement is found only for $k=3$; no numerical improvement over $\ell_k$ was found for $k\in\{4,\ldots,20\}$. The dimension-free upper bound $u=1$ (cyan line) is derived in Theorem~\ref{thm:SA-M2-c1}. Thus, $c_k^\ast\in[s_k, 1]$ for each dimension shown.}
\label{fig:thresholds}
\end{figure}

Among the dimensions investigated, an improvement was found only for $k=3$. Specifically, the parameter vector 
\[
\ba\approx(2.5959305, 0.2602556, 3.2801387),
\]
gives $\Pi(\ba)\approx-8.646209 \cdot10^{-7}$, therefore
\[
c^\ast_3\geq s_3=\min_{i\in[3]}\alpha_i\approx0.2602556>0.2111565\approx\ell_3,
\]
where all values are rounded to seven decimal places. Figure~\ref{fig:thresholds} presents the combined results of the mathematical and computational analysis conducted in this study. Of course, unsuccessful stochastic searches do not prove the absence of counterexamples. The analysis was executed in R version 4.6.1 \citep{R2026}, using the packages \texttt{GenSA} for stochastic optimization \citep{Xiang2013}, \texttt{parallel} for parallel computation \citep{R2026}, and \texttt{ggplot2} for visualization \citep{Wickham2016}. 

\section{Concluding remarks}\label{sec:conclusion}
\noindent
In this study, we focused on three explicit Dirichlet shape parameter estimators and compared them in terms of asymptotic efficiency and componentwise variance. The new moment estimator M2 we propose coincides with M1 for the beta distribution and has smaller componentwise variances in all dimensions $k\geq3$. The results of \citet{tamae2020score} were generalized in two different directions. For the beta distribution, we showed that SA has strictly smaller asymptotic variances than the moment estimators for all $(\alpha, \beta)\in\mathbb{R}_+^2$, not just for $\alpha=\beta$. For the Dirichlet distribution, we showed that SA is asymptotically more efficient than M2 when all parameters are equal and that the difference between SA's and M2's asymptotic variance-covariance matrices is indefinite in any other case.

Finally, we focused on sufficient conditions under which SA is preferable to M2 in terms of componentwise asymptotic variances. SA improves upon M2 when all parameters exceed $c^\ast_k$, i.e., $\ba\in(c^\ast_k,+\infty)^k$, where $c^\ast_k\leq1$.  We also derived a lower bound $\ell_k\leq c^\ast_k$. This bound was computationally improved only for $k=3$. Hence, it remains an open question whether this bound gives the optimal threshold for $k\geq4$. Subsequent studies could work on this problem and refine the analytical characterization of the parameter region in which SA outperforms M2.

\backmatter

\begin{appendices}

\section*{Appendix}\label{appn}

Lemma \ref{lem:trigamma} provides an upper bound for the trigamma function $\psi_1(x)$, which is utilized in Theorem \ref{thm:SA-M1-beta}. \citet{guo2015sharp} also derived this upper bound for the general case of the polygamma functions and showed it is sharp. The proof presented here, however, uses elementary arguments.

\begin{lemma}\label{lem:trigamma}
For every $x\in\mathbb{R}_{+}$,
\begin{equation*}
\frac{1}{x}<\psi_1(x)
<
\frac{1}{x+\tfrac12}+\frac{1}{x^2}.
\end{equation*}
\end{lemma}

\begin{proof}
For the lower bound, let $x, t\in\mathbb{R}_{+}$. Then, the function $t\mapsto(x+t)^{-2}$ is strictly decreasing, therefore
\[
\psi_1(x)
=\sum_{n=0}^{\infty}\frac{1}{(x+n)^2}
>\sum_{n=0}^{\infty}\int_n^{n+1}\frac{1}{(x+t)^2}dt
=\int_0^\infty\frac{1}{(x+t)^2}dt
=\frac1x.
\]
For the upper bound, the series representation of the trigamma function and Tonelli's theorem give
\begin{equation*}
\psi_1(x) =
\sum_{n=0}^{\infty}\frac{1}{(x+n)^2} =
\sum_{n=0}^{\infty}
\int_0^{\infty} t e^{-(x+n)t}\,dt =
\int_0^{\infty}
\sum_{n=0}^{\infty} t e^{-(x+n)t}\,dt =
\int_0^{\infty}
\frac{t e^{-xt}}{1-e^{-t}}\,dt.
\end{equation*}
It is enough to prove that
\begin{equation}\label{eq:kernel-bound}
\frac{t}{1-e^{-t}}<e^{-t/2}+t,
\qquad t>0,
\end{equation}
leading to
\begin{equation*}
\psi_1(x)
= \int_0^{\infty}\frac{t e^{-xt}}{1-e^{-t}}\,dt < \int_0^{\infty}(e^{-t/2}+t)e^{-xt}\,dt = \frac{1}{x+\tfrac12}+\frac{1}{x^2}.
\end{equation*}
Relation \eqref{eq:kernel-bound} is equivalent to
\begin{equation*}
e^{-t/2}-e^{-3t/2}-te^{-t}>0.
\end{equation*}
Setting $s=e^{-t/2}\in(0,1)$ and dividing by $s$ yields
\begin{equation*}
\phi(s):=1-s^2+2s\log s>0.
\end{equation*}
It is easy to prove that $\phi$ is strictly decreasing on $(0,1)$ and $\phi(s)>\phi(1)=0$, therefore \eqref{eq:kernel-bound} holds and the result follows.
\end{proof}

\begin{lemma}\label{lem:beta-gamma-bound}
Let $\boldsymbol{\alpha}\in\mathbb{R}^k_{+}$, $k\geq2$ and $d(\boldsymbol{\alpha}):=3\alpha_0b(\boldsymbol{\alpha})-2c(\boldsymbol{\alpha})$. Then,
\begin{equation*}
0\leq\alpha_0d(\boldsymbol{\alpha})<b^2(\boldsymbol{\alpha}).
\end{equation*}
The left equality holds iff $k=2$.
\end{lemma}

\begin{proof}
Let $e_2$ and $e_3$ be the second and third elementary symmetric polynomials in $\alpha_1,\ldots,\alpha_k$, i.e.,
\begin{equation*}
    e_2 = 1_{k\geq2}\sum_{j=1}^k\sum_{\lambda>j}\alpha_j\alpha_\lambda, \quad e_3 = 1_{k\geq3}\sum_{j=1}^k\sum_{\lambda>j}\sum_{m>\lambda}\alpha_j\alpha_\lambda\alpha_m.
\end{equation*}
Then, $b(\boldsymbol{\alpha})=2e_2$ and $d(\boldsymbol{\alpha})=3\alpha_0b(\boldsymbol{\alpha})-2c(\boldsymbol{\alpha})=6e_3\geq0$, with the equality holding for $k=2$. 
For $k\geq3$, Newton's inequality gives
\begin{equation*}
e_2^2
\geq
\frac{3(k-1)}{2(k-2)}\alpha_0e_3
>
\frac{3}{2}\alpha_0e_3.
\end{equation*}
Thus,
\[
b^2(\boldsymbol{\alpha})=4e_2^2>6\alpha_0e_3=\alpha_0\left\{3\alpha_0b(\boldsymbol{\alpha})-2c(\boldsymbol{\alpha})\right\}=\alpha_0d(\boldsymbol{\alpha}).
\]
\end{proof}

\begin{lemma}\label{lemma:p-positive}
For every $\ba\in\mathbb{R}^k_{+}, \, k\geq2$,
\begin{equation*}
p(\ba)>0.
\end{equation*}
For each fixed $\alpha_0>0$,
\begin{equation}\label{eq:p-infimum}
\inf\left\{p(\ba):\ba\in\mathbb{R}^k_{+},\ \sum_{i=1}^k\alpha_i=\alpha_0\right\}=0.
\end{equation}
\end{lemma}

\begin{proof}

Suppose $0<\alpha_0\leq1$. Let $e_2,e_3$ be the second and third elementary symmetric polynomials in $\alpha_1,\ldots,\alpha_k$. Direct substitution of $b(\boldsymbol{\alpha})=2e_2$ and $c(\boldsymbol{\alpha})=3\alpha_0e_2-3e_3$ in \eqref{eq:p-definition} gives
\begin{equation*}
p(\ba)
=
4(2\alpha_0^2+\alpha_0+1)e_2^2
+4(\alpha_0-1)^2(\alpha_0+1)e_2
-12(\alpha_0-1)(\alpha_0+1)^2e_3.
\end{equation*}
The first term on the right-hand side is strictly positive, and the other two are non-negative, hence $p(\ba)>0$ for $\alpha_0\in(0, 1]$. Suppose $\alpha_0>1$. By Lemma \ref{lem:beta-gamma-bound}:
\begin{equation}\label{eq:ca-inequality}
    2c(\boldsymbol{\alpha})> 3\alpha_0b(\boldsymbol{\alpha})- \frac{b^2(\boldsymbol{\alpha})}{\alpha_0}.
\end{equation}
Substituting \eqref{eq:ca-inequality} into \eqref{eq:p-definition} gives
\[
p(\ba)>
\frac{-\alpha_0^2+3\alpha_0+2}{\alpha_0}b^2(\ba)
+2(\alpha_0+1)(\alpha_0-1)^2b(\ba).
\]
If $1<\alpha_0\le\left(3+\sqrt{17}\right)/2$, the first term on the right-hand side is non-negative and the second is strictly positive. Hence $p(\ba)>0$.

If $\alpha_0>\left(3+\sqrt{17}\right)/2$, the coefficient of $b^2(\ba)$
is negative. Therefore, using $0<b(\ba)<\alpha_0^2$ results in
\[
p(\ba)>
\left\{
\alpha_0(-\alpha_0^2+3\alpha_0+2)
+2(\alpha_0+1)(\alpha_0-1)^2
\right\}b(\ba)=
(\alpha_0^3+\alpha_0^2+2)b(\ba)>0.
\]
Thus $p(\ba)>0$ for every $\alpha_0>1$.
For the infimum, take
\begin{equation*}
\ba_{\varepsilon}
=
\left(\alpha_0-(k-1)\varepsilon,\varepsilon,\ldots,\varepsilon\right),
\qquad 0<\varepsilon<\frac{\alpha_0}{k-1}.
\end{equation*}
Then
\[
\lim_{\varepsilon\rightarrow0^+}b(\ba_{\varepsilon})=\lim_{\varepsilon\rightarrow0^+}c(\ba_{\varepsilon})=0,
\]
so by the definition of $p(\ba)$ in \eqref{eq:p-definition},
\[
    \lim_{\varepsilon\rightarrow0^+}p(\ba_{\varepsilon})
    =0,
\]
completing the proof.
\end{proof}

\begin{lemma}\label{trigamma-inequality-lemma}
For every $x\ge1$,
\[
\psi_1(x)<\frac1x+\frac{2}{3x^2}.
\]
\end{lemma}

\begin{proof}
From Lemma \ref{lem:trigamma},
\[
\psi_1(x)<\frac{1}{x+\tfrac12}+\frac{1}{x^2}=\frac1x+\frac{2}{3x^2}-\frac{x-1}{3x^2(2x+1)}\le\frac1x+\frac{2}{3x^2},
\]
where the last inequality follows from $x\ge1$.
\end{proof}

\begin{lemma}\label{lem:Dirichlet-eigenvalue-sums}
Let $k\ge3$, $\ba\in[1,+\infty)^k$, and
$I:=\sum_{i=1}^k\alpha_i^{-1}$. Define
\[
D=\operatorname{diag}(\alpha_1,\ldots,\alpha_k),
\qquad
M=D-\frac{\ba\ba^{\mathsf T}}{\alpha_0}.
\]
Then, $M$ is positive semidefinite of rank $k-1$ and its
nonzero eigenvalues $\lambda_1,\ldots,\lambda_{k-1}$
satisfy $\lambda_j\ge1$ and
\[
\sum_{j=1}^{k-1}\lambda_j=\frac{b(\ba)}{\alpha_0},
\qquad
\sum_{j=1}^{k-1}\lambda_j^2
=\frac{b^2(\ba)}{\alpha_0^2}-3b(\ba)+\frac{2c(\ba)}{\alpha_0},
\qquad
\sum_{j=1}^{k-1}\lambda_j^{-1}=\frac{k-1}{k}I.
\]
\end{lemma}

\begin{proof}
Let
\[
    v=\left(\sqrt{\alpha_1},\ldots,\sqrt{\alpha_k}\right)^{\mathsf T},
\]
and choose a matrix $H$ whose columns form an orthonormal
basis of $v^\perp$. Since $v^{\mathsf T}v=\alpha_0$,
the orthogonal projection onto $v^\perp$ gives
\[
H^{\mathsf T}H=I_{k-1},
\qquad
HH^{\mathsf T}=I_k-\frac{vv^{\mathsf T}}{\alpha_0}.
\]
Using $D^{1/2}v=\ba$ yields
\[
M=D-\frac{\ba\ba^{\mathsf T}}{\alpha_0}=D^{1/2}HH^{\mathsf T}D^{1/2}
=(D^{1/2}H)(D^{1/2}H)^{\mathsf T}.
\]
Set $A=D^{1/2}H$, so that $M=AA^{\mathsf T}$ and
$A^{\mathsf T}A=H^{\mathsf T}DH$. The form of $M$ implies that it is positive semidefinite with
\[
\operatorname{rank}(M)=\operatorname{rank}(A)=k-1.
\]
Furthermore, $M$ and $H^{\mathsf T}DH$ have the same nonzero
eigenvalues, including multiplicities. Since $A$ has
rank $k-1$, the nonzero eigenvalues
$\lambda_1,\ldots,\lambda_{k-1}$ of $M$ are precisely
the eigenvalues of $H^{\mathsf T}DH$. From the assumption $\ba\in[1,+\infty)^k$,
\[
D\succeq I_k \Longrightarrow H^{\mathsf T}DH\succeq H^{\mathsf T}H=I_{k-1} \Longrightarrow \lambda_j\ge1, \quad j\in[k-1].
\]
Write $s_2=\sum_{i=1}^k\alpha_i^2=\alpha_0^2-b(\ba)$ and
$s_3=\sum_{i=1}^k\alpha_i^3=\alpha_0^3-c(\ba)$. Taking traces gives
\[
\sum_{j=1}^{k-1}\lambda_j
=\operatorname{tr}(M)
=\alpha_0-\frac{s_2}{\alpha_0}
=\frac{b(\ba)}{\alpha_0}.
\]
Similarly,
$\operatorname{tr}\left(D\ba\ba^{\mathsf T}\right)=s_3$ and
$\operatorname{tr}\left\{(\ba\ba^{\mathsf T})^2\right\}=s_2^2$ imply
\[
\sum_{j=1}^{k-1}\lambda_j^2
=\operatorname{tr}\left(M^2\right)
=s_2-\frac{2s_3}{\alpha_0}+\frac{s_2^2}{\alpha_0^2}
=\frac{b^2(\ba)}{\alpha_0^2}-3b(\ba)+\frac{2c(\ba)}{\alpha_0}.
\]
For every nonempty $J\subseteq[k]$, let $M_{J,J}$ denote
the principal submatrix of $M$ indexed by $J$, and write
\[
D_J=\operatorname{diag}(\alpha_i:i\in J),
\qquad
\ba_J=(\alpha_i)_{i\in J}.
\]
Here $\det(D_J)=\prod_{i\in J}\alpha_i$ and
$\ba_J^{\mathsf T}D_J^{-1}\ba_J=\sum_{i\in J}\alpha_i$. Since $D_J$ is invertible, the matrix determinant lemma gives
\begin{equation}\label{eq:detMjj}
\det M_{J,J}=\det(D_J)\left(1-\frac{\ba_J^{\mathsf T}D_J^{-1}\ba_J}{\alpha_0}\right)=\frac1{\alpha_0}\left(\prod_{i\in J}\alpha_i\right)\sum_{\ell\notin J}\alpha_\ell.
\end{equation}
For $1\le r\le n$, define the elementary symmetric polynomial
\[
e_r(x_1,\ldots,x_n)
:=\sum_{1\le i_1<\cdots<i_r\le n}x_{i_1}\cdots x_{i_r}.
\]
These polynomials are related to the principal minors through $\det(I_k+tM)$. Expanding the determinant by
multilinearity in its columns gives
\[
\det(I_k+tM)
=1+\sum_{r=1}^k t^r
 \sum_{\substack{J\subseteq[k]\\|J|=r}}\det M_{J,J}.
\]
On the other hand, since the
eigenvalues of $M$ are $\lambda_1,\ldots,\lambda_{k-1},0$,
\[
\det(I_k+tM)
=\prod_{j=1}^{k-1}(1+t\lambda_j)
=1+\sum_{r=1}^{k-1}
 e_r(\lambda_1,\ldots,\lambda_{k-1})t^r.
\]
Comparing coefficients of $t^r$ therefore yields
\[
e_r(\lambda_1,\ldots,\lambda_{k-1})
=\sum_{\substack{J\subseteq[k]\\|J|=r}}\det M_{J,J},
\qquad 1\le r\le k-1.
\]
Substitution of \eqref{eq:detMjj} gives
\[
e_r(\lambda_1,\ldots,\lambda_{k-1})
=\frac1{\alpha_0}
 \sum_{\substack{J\subseteq[k]\\|J|=r}}
 \sum_{\ell\notin J}
 \left(\prod_{i\in J}\alpha_i\right)\alpha_\ell=\frac{r+1}{\alpha_0}
 \sum_{\substack{S\subseteq[k]\\|S|=r+1}}
 \prod_{i\in S}\alpha_i
=\frac{r+1}{\alpha_0}e_{r+1}(\ba),
\]
since each set $S\subseteq[k]$ of size $r+1$ occurs exactly $r+1$ times in the double sum, corresponding to the choices $\ell\in S$ and $J=S\setminus\{\ell\}$. Finally, 
\[
\sum_{j=1}^{k-1}\lambda_j^{-1}
=\frac{e_{k-2}(\lambda_1,\ldots,\lambda_{k-1})}
       {e_{k-1}(\lambda_1,\ldots,\lambda_{k-1})}
=\frac{(k-1)e_{k-1}(\ba)/\alpha_0}
       {k e_k(\ba)/\alpha_0}
=\frac{k-1}{k}I,
\]
concluding the proof.
\end{proof}

\begin{lemma}\label{lem:Q-mixed}
Let $k\ge3$, $\ba\in[1,+\infty)^k$, and
$I:=\sum_{i=1}^k\alpha_i^{-1}$. Then
\begin{equation}\label{eq:Q-mixed}
\begin{aligned}
&4(k-1)^2\alpha_0c(\ba)(\alpha_0+1)^2(\alpha_0-1)
-\frac23\alpha_0I b^2(\ba)(\alpha_0+2)(\alpha_0+3)\\
&\quad>
6(k-1)^2\alpha_0^2b(\ba)(\alpha_0+1)^2(\alpha_0-1)
-2(k-1)(k-2)b^2(\ba)(\alpha_0+1)^2(\alpha_0-1)\\
&\qquad\quad
-\frac{2k(k-1)}3\alpha_0^2b(\ba)(\alpha_0+2)(\alpha_0+3)
-26(k-1)b^2(\ba).
\end{aligned}
\end{equation}
\end{lemma}

\begin{proof}
Let $\lambda_1,\ldots,\lambda_{k-1}$ be the nonzero
eigenvalues in Lemma~\ref{lem:Dirichlet-eigenvalue-sums},
and define their mean and variance by
\[
\mu=\frac1{k-1}\sum_{j=1}^{k-1}\lambda_j,
\qquad
V=\frac1{k-1}\sum_{j=1}^{k-1}(\lambda_j-\mu)^2.
\]
The eigenvalue identities of Lemma~\ref{lem:Dirichlet-eigenvalue-sums} give
\[
\mu=\frac{b(\ba)}{(k-1)\alpha_0},
\qquad
\frac{V}{\mu^2}
=k-2+\frac{(k-1)\alpha_0
\bigl(2c(\ba)-3\alpha_0b(\ba)\bigr)}{b^2(\ba)}.
\]
Since $k\ge3$ and all eigenvalues are positive,
\[
0\le V
=\frac1{k-1}\sum_{j=1}^{k-1}\lambda_j^2-\mu^2
<\frac1{k-1}\left(\sum_{j=1}^{k-1}\lambda_j\right)^2-\mu^2
=(k-2)\mu^2.
\]
Moreover, the reciprocal-sum identity and $\lambda_j\ge1$
imply
\[
\mu^2\frac{I}{k}-\mu
=\frac1{k-1}\sum_{j=1}^{k-1}
 \frac{(\lambda_j-\mu)^2}{\lambda_j}
\le\frac1{k-1}\sum_{j=1}^{k-1}(\lambda_j-\mu)^2
=V.
\]
Consequently,
\begin{equation}\label{eq:I-bound}
    0\le\frac{V}{\mu^2}<k-2, \qquad I\le\frac{k}{\mu}+\frac{kV}{\mu^2}.
\end{equation}
Let $R(\ba)$ denote the left-hand side of
\eqref{eq:Q-mixed} minus its right-hand side, divided by
$(k-1)b^2(\ba)$. Substituting the expressions for
$\mu$ and $V/\mu^2$ gives
\[
R(\ba)
=26+2(\alpha_0+1)^2(\alpha_0-1)\frac{V}{\mu^2}
+\frac{2\alpha_0(\alpha_0+2)(\alpha_0+3)}{3(k-1)}
 \left(\frac{k}{\mu}-I\right).
\]
The bound on $I$ in \eqref{eq:I-bound} therefore yields
\begin{equation}\label{eq:R-bound}
    R(\ba)\ge 26+\left[2(\alpha_0+1)^2(\alpha_0-1) -\frac{2k\alpha_0(\alpha_0+2)(\alpha_0+3)}{3(k-1)} \right]\frac{V}{\mu^2}.
\end{equation}
If $k=3$, the coefficient in brackets in \eqref{eq:R-bound} equals
\[
-26+(\alpha_0-3)(\alpha_0^2-8).
\]
Since $\alpha_0\ge3$ and $0\le V/\mu^2<1$,
\[
R(\ba)\ge
26\left(1-\frac{V}{\mu^2}\right)
+(\alpha_0-3)(\alpha_0^2-8)\frac{V}{\mu^2}>0.
\]
If $k\ge4$, then $\alpha_0\ge4$ and $k/(k-1)\le4/3$.
Hence the coefficient in brackets in \eqref{eq:R-bound} is bounded below by
\[
2(\alpha_0+1)^2(\alpha_0-1)
-\frac89\alpha_0(\alpha_0+2)(\alpha_0+3)
=\frac{10(\alpha_0-4)^3+98(\alpha_0-4)^2
+238(\alpha_0-4)+6}{9}>0.
\]
Since $V/\mu^2\ge0$, this gives $R(\ba)\ge26>0$.
Thus \eqref{eq:Q-mixed} holds in both cases, and the result follows.
\end{proof}

\begin{lemma}\label{lem:x-trigamma-decreasing}
The function $x\mapsto x\psi_1(x)$ is strictly decreasing on $(0,+\infty)$, with
\[
\lim_{x\rightarrow0^+}x\psi_1(x)=+\infty,
\qquad
\lim_{x\to+\infty}x\psi_1(x)=1.
\]
\end{lemma}

\begin{proof}
The integral representation of the trigamma function and the substitution $u=xt$ give
\begin{equation}\label{eq:x-trigamma-integral}
x\psi_1(x)
=x\int_0^\infty\frac{t e^{-xt}}{1-e^{-t}}\,dt
=\int_0^\infty e^{-u}\frac{u/x}{1-e^{-u/x}}\,du.
\end{equation}
For $y>0$,
\[
\frac{d}{dy}\frac{y}{1-e^{-y}}
=\frac{e^{-y}(e^y-1-y)}{(1-e^{-y})^2}>0.
\]
Thus the integrand in \eqref{eq:x-trigamma-integral} is strictly decreasing in $x$ for every $u>0$, proving strict monotonicity.

The series representation gives
\[
x\psi_1(x)=x\sum_{n=0}^\infty\frac1{(x+n)^2}>\frac1x,
\]
and hence
\[
\lim_{x\rightarrow0^+}x\psi_1(x)=+\infty.
\]
Finally, for $y>0$,
\[
1<\frac{y}{1-e^{-y}}
=y+\frac{y}{e^y-1}<1+y.
\]
Applying these bounds in \eqref{eq:x-trigamma-integral} yields
\[
1<x\psi_1(x)<\int_0^\infty e^{-u}\left(1+\frac ux\right)\,du
=1+\frac1x,
\]
therefore, $\lim_{x\to+\infty}x\psi_1(x)=1$, concluding the proof.
\end{proof}

\end{appendices}

\bibliography{refs}

\end{document}